\documentclass[12pt]{article}
\usepackage[mathscr]{eucal}
\usepackage{amssymb,amsmath,amscd}
\usepackage{eufrak}
\font\cmc=cmcsc10  scaled \magstep2

\newcommand\F{\mathbb{F}}

\newcommand\mscrP{\mathscr{P}}

\newcommand\vk{\vskip}

\newcommand\al{\alpha}
\newcommand\be{\beta}

\newcommand\und{\underset}
\newcommand\no{\noindent}
\newcommand\unl{\underline}

\begin{document}
\vbox to .5truecm{}

\begin{center}
\cmc
A note on Je{\'s}manowicz' conjecture for non-primitive Pythagorean triples
\end{center}
\vspace{0.5cm}

\begin{center}
Van Thien Nguyen,\ Viet Kh. Nguyen,\ Pham Hung Quy
\end{center}

\begin{abstract} {Let $(a, b, c)$ be a primitive Pythagorean triple parameterized as $a=u^2-v^2,\ b=2uv,\ c=u^2+v^2$,\ where $u>v>0$ are co-prime and 
not of the same parity. In 1956, L. Je{\'s}manowicz conjectured that for any positive integer $n$, the Diophantine equation $(an)^x+(bn)^y=(cn)^z$ has 
only the positive integer solution $(x,y,z)=(2,2,2)$. In this connection we call a positive integer solution $(x,y,z)\ne (2,2,2)$ with $n>1$ exceptional. In 1999 M.-H. Le 
gave necessary conditions for the existence of exceptional solutions which were refined recently by H. Yang and R.-Q. Fu. In this paper we give a unified simple proof 
of the theorem of Le-Yang-Fu. Next we give necessary conditions for the existence of exceptional solutions in the case $v=2,\ u$ is an odd prime. As an application we show 
the truth of the Je{\'s}manowicz conjecture for all prime values $u < 100$.}

\end{abstract}



\begin{center}
{\bf 1.\ Introduction}
\end{center}

Let $(a, b, c)$ be a primitive Pythagorean triple. Clearly for such a triple with $2\mid b$ one has the following parameterization
  $$a=u^2-v^2,\ b=2uv,\ c=u^2+v^2$$
with
$$u>v>0,\ \mathrm{gcd}(u,v)=1,\ u+v\equiv 1\ (\mathrm{mod}\ 2)\leqno{(1.1)}$$

In 1956 L. Je{\'s}manowicz (\cite{Jes56}) made the following
\vk.3cm
{\bf Conjecture\ 1.1.}\quad {\it For any positive integer $n$, the Diophantine equation
  $$(an)^x+(bn)^y=(cn)^z\leqno{(1.2)}$$
has only the positive integer solution $(x,y,z)=(2,2,2)$.}
\vk.2cm
The primitive case of the conjecture ($n=1$) was investigated thoroughly. Although the conjecture is still open, many special cases are shown to be true. We refer to a recent survey \cite{LeSS18} for a detailed account.
\vk.2cm
Much less known about the non-primitive  case ($n>1$). A positive integer solution $(x, y, z, n)$ of $(1.2)$ is called exceptional if $(x, y, z)\ne (2, 2, 2)$ and $n > 1$. 
For a positive integer $t$ let $\mscrP(t)$ denote the set of distinct prime factors of $t$, and $P(t) -$ their product.
The first known result in this direction was obtained in 1998 by M.-J. Deng and G.L. Cohen (\cite{DeCo98}), namely if $u=v+1$,\ $a$ is a prime power, and either $P(b)\mid n$, 
or $P(n)\nmid b$, then $(1.2)$ has only positive integer solution $(x,y,z)=(2,2,2)$. 
In 1999, M.-H. Le gave necessary conditions for $(1.2)$ to have exceptional solutions.
\vk.3cm
{\bf Theorem\ 1.2.}\ (\cite{Le99})\quad {\it If $(x,y,z,n)$ is an exceptional solution of $(1.2)$, then one of the following three conditions is satisfied:
\vk.1cm
$(i)$\ $\max\{ x,y\}>\min\{ x,y\}>z$,\ $\mscrP(n)\subsetneqq \mscrP(c);$
\vk.1cm
$(ii)$\ $x>z>y$,\ $\mscrP(n)\subset \mscrP(b);$
\vk.1cm
$(iii)$\ $y>z>x$,\ $\mscrP(n)\subset \mscrP(a).$}
\vk.2cm
However, as noted in \cite{YF15} by H. Yang and R.-Q. Fu, the case $x=y>z$ is not completely handled by the arguments used in \cite{Le99}. Furthermore they completed the unhandled case (\cite{YF15}, Theorem 1) based on 
a powerful result of Zsigmondy (\cite{Zs92}, {\it cf.} \cite{BV04}, \cite{Car13}).
In fact one can give a unified simple proof of Theorem of Le-Yang-Fu (Theorem 1.2) by using a weaker version of the Zsigmondy theorem as stated in Lemma 3 of \cite{DeCo98} ({\it cf.} \S 2 below).
\vk.2cm
Since \cite{DeCo98}, \cite{Le99} many works intensively investigated the first interesting family of primitive triples:
 $$v=1,\ u=2^k,\ k=1, 2,\ldots\leqno{(1.3)}$$ 

Most recently X.-W. Zhang and W.-P. Zhang (\cite{ZZ14}), and T. Miyazaki (\cite{Miy15}) independently proved Conjecture 1.1 for the (infinite) family $(1.3)$.
\vk.2cm
It is natural to treat the next interesting case:\ $v=2,\ u$ is an {\it odd prime} which was known recently for few values $u$:\ $u=3$ (\cite{DeCo98}),\ $u=5$ -- by Z. Cheng, 
C.-F. Sun and X.-N. Du,\ $u=7$ -- by C.-F. Sun, Z. Cheng, and by G. Tang,\ $u=11$ -- by W.-Y. Lu, L. Gao and H.-F. Hao ({\it cf.} \cite{LeSS18} for references). 
Let's formulate our main results. We rewrite $(1.2)$ as
  $$[(u^2-4)n]^x+(4un)^y=[(u^2+4)n]^z\leqno{(1.4)}$$

An arithmetical argument (given in Lemma 3.5 below) shows that $u^2-4$ admits a proper decomposition $u^2-4=u_1 u_2,\ \gcd(u_1,u_2)=1$, so that there are three possibilities to consider: $u_1\equiv \pm 1, 5\ (\mathrm{mod}\ 8)$. 
\vk.3cm
{\bf Theorem\ 1.3.}\quad {\it If $(x,y,z,n)$ is an exceptional solution of $(1.4)$ and $u_1\equiv \pm 1\ (\mathrm{mod}\ 8)$, then $y$ is even.}
\vk.2cm
In view of Theorem 1.3 the possibility $u_1\equiv -1\ (\mathrm{mod}\ 8)$ is eliminated, because in this case $x, y, z$ are even, which is in general impossible by an auxiliary argument (Lemma 3.6 below).
\vk.2cm
Let $\nu_q(t)$, for a prime $q$,\ denote the exponent of $q$ in the prime factorization of $t$, and let $\Big( \dfrac{}{m}\Big)$ denote the Jacobi quadratic residue symbol. 
\vk.3cm
{\bf Theorem\ 1.4.}\quad {\it If $(x,y,z,n)$ is an exceptional solution of $(1.4)$, then one of the following cases is satisfied
\vk.2cm
$(1)$\ $\nu_2(u_1-1)=3$$:$\ $(\nu_2(x), \nu_2(y), \nu_2(z))=(0, \geq 2, 1);$\ $u_1$ admits a proper decomposition $u_1=t_1 t_2,\ \gcd(t_1,t_2)=1$ and $t_1, t_2\equiv 5\ (\mathrm{mod}\ 8)$ satisfying certain 
special Diophantine equations$;$
\vk.2cm
$(2)$\ $u_1\equiv 5\ (\mathrm{mod}\ 8)$,\ $u_2=w^{2^s}$, where $s=\nu_2(z-x)-\nu_2(x)$ and either of the following
\vk.1cm
\quad $(2.1)$\ $w\equiv \pm 3\ (\mathrm{mod}\ 8)$$:$\ $(\nu_2(x), \nu_2(y), \nu_2(z))=(0, \geq 1, 0);$\ $u\equiv 1\ (\mathrm{mod}\ 4);$\ $\Big( \dfrac{u_1}{p}\Big)=\Big( \dfrac{w}{p}\Big),\ \forall\ p\mid (u^2+4)$ and 
$\Big( \dfrac{w}{p}\Big)=\Big( \dfrac{u^2+4}{p}\Big),\ \forall\ p\mid u_1;$
\vk.2cm
\quad $(2.2)$\ $w\equiv \pm 1\ (\mathrm{mod}\ 8)$$:$\ $(\nu_2(x), \nu_2(y), \nu_2(z))=(\be, 0, \be),\ \be\geq 1;$\ $u\equiv\pm 3\ (\mathrm{mod}\ 8);$\ $\Big( \dfrac{w}{p}\Big)=1,\ \forall\ p\mid (u^2+4)$ and 
$\Big( \dfrac{w}{p}\Big)=\Big( \dfrac{u}{p}\Big),\ \forall\ p\mid u_1$. Moreover, if $u\equiv 3\ (\mathrm{mod}\ 8)$,\ then $w$ can not be a square.}

\vk.3cm
{\bf Corollary\ 1.5.}\quad {\it Conjecture 1.1 is true for\ $v=2$,\ $u$ -- an odd prime $< 100$.}
\vk.2cm
Let's explain the ideas in proving our main results. As for Theorem 1.3 and Theorem 1.4 we exploit a total analysis of Jacobi quadratic and quartic residues. In the case $u_1\equiv 1\ (\mathrm{mod}\ 8)$\ we have 
a further proper decomposition $u_1=t_1 t_2$, which leads to certain special Diophantine equations. Theorem 1.4 helps us substantially in reducing the 
verification process, as the possibility $u_1\equiv 5\ (\mathrm{mod}\ 8)$ occurs quite sparsely. We demonstrate this for $u < 100$ in proving Corollary 1.5.
\vk.2cm
The paper is organized as follows. In \S 2 we give a unified simple proof of Theorem 1.2. Section 3 provides some reduction of the problem and preliminary results. 
Theorem 1.3 will be proved  in \S 4. The case $u_1\equiv 5\ (\mathrm{mod}\ 8)$ and Theorem 1.4 will be treated in \S 5. The verification for $u < 100$ in 
Corollary 1.5 will be given in the final section.
\vk.3cm
{\bf Acknowledgement.}\quad The authors would like to thank the referee for many valuable comments and suggestions greatly improving the content of the paper. 
\vk.5cm

\begin{center}
{\bf 2.\ A Simple Proof of Theorem 1.2}
\end{center}

We shall use the following weaker version of Zsigmondy's theorem.
\vk.3cm
{\bf Lemma\ 2.1} ({\it cf.} \cite{DeCo98}, Lemma 3).\quad {\it For $X>Y>0$ co-prime integers,
\vk.1cm
$(1)$\ if $q$ is a prime, then 
  $$\gcd\big( X-Y,\dfrac{X^q-Y^q}{X-Y}\big)=1,\ {\text{\rm or}}\ q; $$

$(2)$\ if $q$ is an odd prime, then 
  $$\gcd\big( X+Y,\dfrac{X^q+Y^q}{X+Y}\big)=1,\ {\text{\rm or}}\ q. $$
}

Part $(2)$ is Lemma 3 of \cite{DeCo98}. As for part $(1)$ one argues similarly:\ if $\ell^r$ is a common prime power divisor of $X-Y$ and $(X^q-Y^q)/(X-Y)$. Clearly 
    $$\dfrac{X^q-Y^q}{X-Y}\equiv 0\ (\mathrm{mod}\ \ell^r)\leqno{(2.1)}$$

On the other hand from the fact that $X\equiv Y\ (\mathrm{mod}\ \ell^r)$ it follows
  $$\dfrac{X^q-Y^q}{X-Y}=X^{q-1}+X^{q-2} Y+\cdots +X Y^{q-2}+Y^{q-1}\equiv q Y^{q-1}\ (\mathrm{mod}\ \ell^r)\leqno{(2.2)}$$

Since $\ell\nmid Y$,\ $(2.1)-(2.2)$ imply that $\ell=q$, and $r=1$. \hfill\ensuremath{\square}
\vk.3cm
{\bf Remark\ 2.2.}\quad Part $(1)$ is a special case of Theorem IV in \cite{BV04}.
\vk.3cm
{\bf Lemma\ 2.3.}\quad {\it For a prime divisor $q$ of $(X-Y)$ and positive integer $\be$
  $$\nu_q(X^{q^\be}-Y^{q^\be})=\be+\nu_q(X-Y)\leqno{(2.3)}$$
}
\vk.2cm
{\it Proof.}\quad Applying part $(1)$ of Lemma 2.1 $\be$ times one has
  $$\gcd\big( X^{q^{\be-1}}-Y^{q^{\be-1}},\dfrac{X^{q^\be}-Y^{q^\be}}{X^{q^{\be-1}}-Y^{q^{\be-1}}}\big)=q; $$
  $$\cdots$$
  $$\gcd\big( X-Y,\dfrac{X^q-Y^q}{X-Y}\big)=q.$$
Hence the formula $(2.3).$ \hfill\ensuremath{\square}
\vk.3cm
In view of Lemma 2 of \cite{DeCo98} there are no exceptional solutions with $z\geq\max\{x,y\}$, so as in \cite{Le99} we have to eliminate the following three cases:

\begin{center} 
{\bf (I)}\ $x>y=z;$\qquad {\bf (II)}\ $y>x=z;$\qquad {\bf (III)}\ $x=y>z.$
\end{center}

\vk.3cm
{\bf (I)}\ \unl{$x>y=z$}:\ dividing both sides of $(1.2)$ by $n^y$ one gets
  $$a^x n^{x-y}=c^y-b^y\leqno{(2.4)}$$

By considering $\mathrm{mod}\ c+b$, and taking into account $(c+b) (c-b)=a^2$, one sees that $y$ must be even, say $y=2 y_1$. 
Now put $X=c^2,\ Y=b^2$, so $X\equiv Y (\mathrm{mod}\ a^2),\ \gcd(Y,a)=1$. Taking $\mathrm{mod}\ a$ and in view of $(2.4)$ 
  $$0\equiv \dfrac{X^{y_1}-Y^{y_1}}{X-Y}=X^{y_1-1}+X^{y_1-2} Y+\cdots +X Y^{y_1-2}+Y^{y_1-1}\equiv y_1 Y^{y_1-1}\ (\mathrm{mod}\ a)$$
one concludes that $a\mid y_1$.
\vk.2cm
For any $q\in\mscrP(a)$ let $\be=\nu_q(y_1)$, so that $y_1=q^\be y_2$ with $q\nmid y_2$. Putting $U=X^{q^\be},\ V=Y^{q^\be}$ for short, we have
  $$X^{y_1}-Y^{y_1}=(U-V) (U^{y_2-1}+U^{y_2-2} V+\cdots +U V^{y_2-2}+V^{y_2-1})\leqno{(2.5)}$$
and
  $$U^{y_2-1}+U^{y_2-2} V+\cdots +U V^{y_2-2}+V^{y_2-1}\equiv y_2 V^{y_2-1}\not\equiv 0\ (\mathrm{mod}\ q) \leqno{(2.6)}$$

Lemma 2.3 and $(2.5), (2.6)$ imply that 
 $$\nu_q(X^{y_1}-Y^{y_1})=\nu_q(U-V)=\be+2\nu_q(a) \leqno{(2.7)}$$
In view of $(2.4)$ the equality $(2.7)$ means that $a^{x-2}\mid y_1$ in contradiction with 
$y_1=y/2<a^{x-2}$ as $x>y,\ a>1$.
\vk.3cm
{\bf (II)}\ \unl{$y>x=z$}:\ Similarly dividing both sides of $(1.2)$ by $n^z$ one gets
  $$b^y n^{y-x}=c^x-a^x\leqno{(2.8)}$$

Arguing as above with $\mathrm{mod}\ c+a$, one sees that $x$ must be even, say $x=2 x_1$. Put $X=c^2,\ Y=a^2$. Considering $\mathrm{mod}\ b$ and from  
$(2.8)$ it follows that $b\mid x_1$. So $\nu_q(X^{x_1}-Y^{x_1})=\nu_q(x_1)+2\nu_q(b)$  for any $q\in\mscrP(b)$, therefore $b^{y-2}\mid x_1$ in contradiction with 
$x_1=x/2<b^{y-2}$ as $y>x,\ b>1$.
\vk.3cm
{\bf (III)}\ \unl{$x=y>z$}:\ dividing both sides of $(1.2)$ by $n^z$ one gets
  $$(a^x+b^x) n^{x-z}=c^z\leqno{(2.9)}$$

First we claim that $x$ must be {\it even}. Indeed, if $x$ is odd, then from $(2.9)$ it follows that there is an odd prime $q\in \mscrP(a+b)\cap\mscrP(c)$, so $q\in\mscrP(ab)$, as 
$c^2=a^2+b^2$. A contradiction with $\gcd(a,b)=1$.
\vk.2cm
Writing now $x=2 x_1$ one sees that $x_1$ must be {\it odd}. Since otherwise for an odd prime $q\in \mscrP(a^x+b^x)\cap\mscrP(c)$ taking $\mathrm{mod}\ q$ and by $(2.9)$
  $$0\equiv a^x+b^x=a^{2 x_1}+(c^2-a^2)^{x_1}\equiv 2 a^{2 x_1} \ (\mathrm{mod}\ q)$$
one gets a contradiction with $\gcd(a,c)=1$.
\vk.2cm
Now from $(2.9)$ we see that
  $$\dfrac{(a^2)^{ x_1}+(b^2)^{x_1}}{a^2+b^2}=\dfrac{c^{z-2}}{n^{x-z}}>1 \leqno{(2.10)}$$
as $x>z\geq 2$. So there is an odd prime $q\in\mscrP(c)$ dividing $((a^2)^{ x_1}+(b^2)^{x_1})/(a^2+b^2)$. Considering $\mathrm{mod}\ q$ and taking into account 
$a^2\equiv -b^2\ \mathrm{mod}\ q,\ q\nmid a$ one has
  $$0\equiv \dfrac{(a^2)^{ x_1}+(b^2)^{x_1}}{a^2+b^2}=(a^2)^{x_1-1}-(a^2)^{x_1-2} b^2+\cdots-a^2 (b^2)^{x_1-2}+(b^2)^{x_1-1}\equiv x_1 a^{2 x_1-2}\ (\mathrm{mod}\ q)$$
Hence $q\mid x_1$, and so $((a^2)^q+(b^2)^q)\mid ((a^2)^{ x_1}+(b^2)^{x_1})$. Applying part $(1)$ of Lemma 2.1 we get
  $$\gcd\big( a^2+b^2,\dfrac{(a^2)^q+(b^2)^q}{a^2+b^2}\big)=q\leqno{(2.11)}$$

On the other hand from $(2.10)$  one knows that $((a^2)^q+(b^2)^q)/(a^2+b^2)$ is a product of primes in $\mscrP(c)$.
It is easy to see that $((a^2)^q+(b^2)^q)/(a^2+b^2)>q$. So either $\nu_q\big(((a^2)^q+(b^2)^q)/(a^2+b^2)\big)\geq 2$ and $\nu_q(a^2+b^2)\geq 2$, or both of them 
must have another common prime factor in $\mscrP(c)$, a contradiction with $(2.11)$.
 
\vk.5cm

\begin{center}
{\bf 3.\ Preliminary reduction}
\end{center}

We need some reduction of the problem. The following result is due to N. Terai.
\vk.3cm
{\bf Lemma\ 3.1} (\cite{Te14}).\quad {\it Conjecture 1.1 is true for $n=1,\ v=2$.} 
\vk.2cm
Because of Lemma 3.1 we will assume henceforth $n>1$.
\vk.2cm
M.-J. Deng (\cite{De14}, from the proof of Lemma 2.3), and H. Yang, R.-Q. Fu (\cite{YF15}) showed that we can remove the condition $(i)$ in Theorem 1.2.
\vk.3cm
{\bf Lemma\ 3.2}\quad {\it If $(x, y, z,n)$ is an exceptional solution, then either $x>z>y$, or $y>z>x.$}
\vk.2cm
Note that the proof of Lemma 3.2 relies essentially on the condition $n>1$. It could be interesting to find a proof of this result for the case $n=1$. 
\vk.2cm
Furthermore, in the case when $u$ is an odd prime and $v=2$, H. Yang, R.-Q. Fu \cite{YRF17} succeeded to eliminate the possibility 
$(ii)$ in Theorem 1.2.
\vk.3cm 
{\bf Lemma\ 3.3.}\quad {\it Suppose that $u$ is an odd prime and $v=2$. Then equation $(1.2)$ has no exceptional solutions $(x, y, z,n)$ with $x > z > y$.}

\vk.3cm
{\bf Lemma\ 3.4.}\quad {\it For a positive integer $w$
\vk.1cm
$(1)$\ if $\nu_2(w)\geq 2,$ then $\nu_2[(1+w)^x-1]=\nu_2(w)+\nu_2(x);$
\vk.1cm
$(2)$\ if $\nu_2(w)=1$ and $x$ is odd, then $\nu_2[(1+w)^x-1]=1;$
\vk.1cm
$(3)$\ if $\nu_2(w)=1$ and $x$ is even, then $\nu_2[(1+w)^x-1]=\nu_2(2+w)+\nu_2(x).$ 
\vk.1cm
In particular $\nu_2[(1+w)^x-1]=2+\nu_2(x),$ if $w\equiv 4\ (\mathrm{mod}\ 8);$ or if $w\equiv 2\ (\mathrm{mod}\ 8)$ and $x$ is even.
}
\vk.2cm
{\it Proof.}\quad $(1)$\ The conclusions of Lemma 3.4 are true trivially for $x=1$. Assuming now $x\geq 2$ we have
  $$(1+w)^x-1=w (C^1_x+C^2_x w+\cdots +C^{x-1}_x w^{x-2}+C^x_x w^{x-1})\leqno{(3.1)} $$

Clearly $\nu_2(j)\leq j-1$ for $j=2, \cdots, x,$ and so
  $$\nu_2(C^j_x w^{j-1})=\nu_2\Big( \dfrac{x}{j} C^{j-1}_{x-1} w^{j-1}\Big) \geq \nu_2(x)+j-1>\nu_2(x),$$
as $\nu_2(w)\geq 2$. Hence the conclusion follows from taking $\nu_2(\cdot)$ on both sides of $(3.1).$
\vk.2cm
$(2)$\ Obvious from $(3.1)$, since $C^1_x+C^2_x w+\cdots +C^{x-1}_x w^{x-2}+C^x_x w^{x-1}$ is {\it odd} in this case.
\vk.2cm
$(3)$\ Writing $x=2 x_1$ we have
  $$(1+w)^x-1=[(1+w)^{x_1}-1] [(1+w)^{x_1}+1] \leqno{(3.2)}$$

If $x_1$ is {\it odd}, {\it i.e.} $\nu_2(x)=1$, then $\nu_2[(1+w)^{x_1}-1]=1$ by the part $(2)$ above, and $\nu_2[(1+w)^{x_1}+1]=\nu_2(2+w)$, as 
  $$(1+w)^{x_1}+1=(2+w) [(1+w)^{x_1-1}-(1+w)^{x_1-2}+\cdots -(1+w)+1] $$
and $(1+w)^{x_1-1}-(1+w)^{x_1-2}+\cdots -(1+w)+1$ is {\it odd}.
\vk.2cm
If $x_1$ is {\it even}, then $\nu_2[(1+w)^{x_1}+1]=1$, since
  $$(1+w)^{x_1}+1=2+C^1_{x_1} w+C^2_{x_1} w^2+\cdots +C^{x_1-1}_{x_1} w^{x_1-1}+w^{x_1}. $$
Therefore $\nu_2[(1+w)^x-1]=\nu_2[(1+w)^{x_1}-1]+1$ by $(3.2)$. Now the descending argument yields the conclusion. \hfill\ensuremath{\square}
\vk.3cm
The following claims play a central role in the next sections. 
\vk.3cm
{\bf Lemma\ 3.5.}\quad {\it If $(x,y,z,n)$ is an exceptional solution of $(1.4)$, then $u^2-4$ admits a proper decomposition $u^2-4=u_1 u_2,\ \gcd(u_1,u_2)=1$ and with one of the following conditions satisfied:
\vk.1cm
$(1)$\ $u_1\equiv 1\ (\mathrm{mod}\ 8)$ and\ $\nu_2(z)=\nu_2(u_1-1)+\nu_2(x)-2;$
\vk.1cm
$(2)$\ $u_1\equiv 7\ (\mathrm{mod}\ 8)$,\ $\nu_2(z)=\nu_2(u_1+1)+\nu_2(x)-2$,\ and $\nu_2(x)\geq 1;$
\vk.1cm
$(3)$\ $u_1\equiv 5\ (\mathrm{mod}\ 8)$,\ $u_2$ is a square and and $\nu_2(z)=\nu_2(x).$}
\vk.3cm
{\it Proof.}\quad In view of Lemmas 3.2, 3.3 we may assume the existence of an exceptional solution with $y>z>x$ (the case $(iii)$ of Theorem 1.2). Dividing both sides of $(1.4)$ by 
$n^x$ one gets
  $$(u^2-4)^x=[(u^2+4)^z-(4u)^yn^{y-z}] n^{z-x} \leqno{(3.3)}$$ 

It is easy to see that $\gcd(u^2+4, n)=1$. So $(3.3)$ is equivalent to the following system
  $$\begin{cases}
     (u^2+4)^z-(4u)^y n^{y-z} &=u_1^{x} \\
     n^{z-x} &=u_2^{x}
  \end{cases}\leqno{(3.4)}$$  
with $u^2-4=u_1 u_2,\ \gcd(u_1,u_2)=1.$ The system $(3.4)$ can be rewritten as
  $$(u^2+4)^z-2^{2y} u^y n^{y-z}=u_1^x \leqno{(3.5)}$$
or equivalently
  $$[(u^2+4)^z-1]-(u_1^x-1)=2^{2y} u^y n^{y-z} \leqno{(3.5')}$$
with $k(z-x)=mx$, and $n^m=u_2^k$.
\vk.2cm
Clearly $u_2>1$. Assume now $u_1=1$. As $u^2\equiv 1\ (\mathrm{mod}\ 8)$, by comparing $\nu_2(\cdot)$ both sides of $(3.5)$ and by $(1)$ of Lemma 3.4 we have $\nu_2[(u^2+4)^z-1]=2+\nu_2(z)<2y$. 
So $(3.5')$ is inconsistent. So $u_1>1$ and 
  $$\nu_2(z)=\nu_2(u_1^x-1)-2 \leqno{(3.6)}$$
\vk.2cm
If $u_1\equiv 1\ \mathrm{mod}\ 8$, then by $(1)$ of Lemma 3.4 we get $\nu_2(z)=\nu_2(u_1-1)+\nu_2(x)-2$.
\vk.2cm
If $u_1 \equiv 7\ \mathrm{mod}\ 8$ and $x$ is {\it odd}, then by $(2)$ of Lemma 3.4:\ $\nu_2(u_1^x-1)=1$, impossible by $(3.6)$. Thus $(3.5')$ is inconsistent. 
\vk.2cm
If $u_1\equiv 7\ \mathrm{mod}\ 8$ and $x$ is {\it even}, then by $(3)$ of Lemma 3.4:\ $\nu_2(u_1^x-1)=\nu_2(u_1+1)+\nu_2(x)$. Hence by $(3.6)$ one gets $\nu_2(z)=\nu_2(u_1+1)+\nu_2(x)-2$.
\vk.2cm
For $u_1\equiv 3\ \mathrm{mod}\ 8$, we have $\nu_2(u_1^x-1)=1$, if $x$ is {\it odd} (by $(2)$ of Lemma 3.4), and $\nu_2(u_1^x-1)=2+\nu_2(x)$, if $x$ is {\it even} (by $(3)$ of Lemma 3.4). 
Hence for $(3.5)$ to be consistent one has necessarily $\nu_2(z)=\nu_2(x)$, which implies $\nu_2(z-x)\geq \nu_2(x)+1$. So from the second equation of 
$(3.4)$:\ $n^{z-x}=u_2^x$ it follows that $u_2$ must be a square, hence $u_2\equiv 1\ \mathrm{mod}\ 8$. Thus $u_1 u_2\equiv 3\mod 8$, a contradiction with 
$u_1 u_2=u^2-4\equiv 5\ \mathrm{mod}\ 8$. 
\vk.2cm
Similarly, for $u_1\equiv 5\ \mathrm{mod}\ 8$, by using $(1)$ of Lemma 3.4 we have $\nu_2(u_1^x-1)=2+\nu_2(x)$, and by the same reason $\nu_2(z)=\nu_2(x)$. Hence the system $(3.4)$ is inconsistent,  
if $u_2$ is not a square. \hfill\ensuremath{\square}
\vk.3cm
{\bf Lemma\ 3.6.}\quad {\it In the notations above if $x, y, z$ are even, then $(3.5)$ is inconsistent.}
\vk.2cm
{\it Proof}.\quad In this case we can rewrite $(3.5)$ in the form of Pythagorian equation
  $$\big(u_1^{x/2}\big)^2+\big[2^y u^{y/2} n^{(y-z)/2}\big]^2=\big[(u^2+4)^{z/2}\big]^2.$$

Hence ({\it cf.} $(1.1)$) there are integers $X, Y$, say with $2\mid Y$ such that
  $$(u^2+4)^{z/2}=X^2+Y^2\leqno{(3.7)}$$
  $$2^y u^{y/2} n^{(y-z)/2}=2 X Y\leqno{(3.8)}$$

In view of \cite{ZZ14}, Lemma 2.2, equation $(3.7)$ has solutions
  $$u^2+4=A^2+B^2,\ \ 2\mid B\leqno{(3.9)}$$
  $$\nu_2(Y)=\nu(z/2)+\nu_2(B)\leqno{(3.10)}$$

Since $u^2+4\equiv 5\ \mathrm{mod}\ 8$ it follows from $(3.9)$ that $\nu_2(B)=1$. From $(3.8)$ we have $\nu_2(Y)=y-1$ which together with $(3.10)$ implies
  $$y=\nu_2(z)+1$$
a contradiction with $y>z$. \hfill\ensuremath{\square}
 
\vk.3cm
{\bf Corollary\ 3.7.}\quad {\it In the notations above if $y, z$ are even and $(3.5)$ is consistent, then $x$ is odd and $u_1\equiv 1\ (\mathrm{mod}\ 8)$. Moreover 
$u_1$ admits a proper decomposition $u_1=t_1 t_2$ such that $\gcd(t_1,t_2)=1$ and
  $$t_2^x+ t_1^x=2  (u^2+4)^{z/2}\leqno{(3.11)}$$
  $$t_2^x-t_1^x=2^{y+1} u^{y/2} n^{(y-z)/2}\leqno{(3.12)}$$
  $$\nu_2(t_1^x-1)=\nu_2(t_2^x-1)=\nu_2(u_1^x-1)-1\leqno{(3.13)}$$
}

{\it Proof}.\quad By Lemma 3.6 $x$ is odd. In fact one can rewrite $(3.5)$ as
  $$A\cdot B=u_1^x\ \mathrm{with}\ \gcd(A,B)=1, $$
where
  $$A=(u^2+4)^{z/2}-2^y u^{y/2} n^{(y-z)/2},\  B=(u^2+4)^{z/2}+2^y u^{y/2} n^{(y-z)/2}.$$

Hence
  $$A=t_1^x,\ B=t_2^x\ \mathrm{with}\ u_1=t_1 t_2\ \mathrm{and}\ \gcd(t_1,t_2)=1 \leqno{(3.14)}$$

If $t_1=1$, then by $(1)$ of Lemma 3.4: $\nu_2[(u^2+4)^{z/2}-1]=2+\nu_2(z/2)<y=\nu_2(2^y u^{y/2} n^{(y-z)/2})$. So $A=1$ is impossible. 
\vk.2cm
Now from $(3.14)$ we have two possibilities:
\vk.1cm
\quad $(1)$\ $z/2$ is odd: $t_1\equiv t_2\equiv 5\ (\mathrm{mod}\ 8);$ \qquad $(2)$\ $z/2$ is even: $t_1\equiv t_2\equiv 1\ (\mathrm{mod}\ 8);$
\vk.2cm
\no both of them imply $u_1\equiv 1\ (\mathrm{mod}\ 8)$. 
\vk.2cm
Also $(3.11)-(3.13)$ follow immediately from $(3.14)$. \hfill\ensuremath{\square}
\vk.3cm
{\bf Corollary\ 3.8.}\quad {\it In the situation of Corollary 3.7 we have $t_1, t_2\equiv 5\ (\mathrm{mod}\ 8)$ and $\nu_2(u_1-1)=3.$}
\vk.2cm
{\it Proof}.\quad We will show that $z/2$ must be odd, from which the conclusion immediately follows by the proof above, noting that $\nu_2(u_1-1)=\nu_2(u_1^x-1)=\nu_2(A-1)+1=3$. 
\vk.2cm
Assume on the contrary that $\nu_2(z)\geq 2$. In view of $(3.14)$ one has $x\geq 3$, as $t_1<t_2<u^2-4$. We claim that $x>3$. Indeed, if $x=3$, then $n=u_2^3$ by $(3.4)$, noting that $z=4$ by $B=t_2^x$ of 
$(3.14)$, so $y=6$ as $A=t_1^x>0$. Now from 
the equation $t_1^x=A$ in $(3.14)$ we see that $(t_1,4 u u_2,u^2+4)$ is a primitive solution of 
     $$X^3+Y^3=Z^2\leqno{(3.15)}$$

Euler (\cite{Dick2}, pp. 578--579) indicated a primitive parameterization for the Diophantine equation $(3.15)$ with $3\nmid Z,\ 2\mid Y$ as follows
  $$X=(s-t) (3s-t) (3s^2+t^2),\ \ Y=4 s t (3s^2-3st+t^2) $$
with $s, t$ co-prime,\ $3\nmid t$ and $s\not\equiv t\ (\mathrm{mod}\ 2)$. Hence $8\mid Y$ which shows that $t_1^x=A$ in $(3.14)$ is impossible. 
\vk.2cm
Furthermore, if $x\geq 4$, then by Theorem 1.1 of \cite{BeS04}, $(3.11)$ is again impossible. \hfill\ensuremath{\square}

\vk.5cm

\begin{center}
{\bf 4.\ Proof of Theorem 1.3}
\end{center}

The aim of this section is to show that the case $u_1\equiv 7\ (\mathrm{mod}\ 8)$ in Lemma 3.5 is not realized.
We refer the reader to \cite{IR82} for basic properties of Jacobi quadratic and quartic residue symbols $\Big( \dfrac{}{m}\Big),\ \Big( \dfrac{}{m}\Big)_4$ we shall use in the following lemmas.
\vk.3cm
{\bf Lemma\ 4.1.}\quad {\it For a prime $p\mid (u^2+4)$ one has $p\equiv 1\ (\mathrm{mod}\ 4)$ and $\Big( \dfrac{u}{p}\Big)=1$.}
\vk.2cm
{\it Proof}.\quad Since $u^2\equiv -4\ (\mathrm{mod}\ p)$, so $\Big( \dfrac{-1}{p}\Big)=1$, {\it i.e.} $p\equiv 1\ (\mathrm{mod}\ 4)$. Furthermore we include the following simple argument due to the referee instead of ours in 
the original version:
  $$\Big( \dfrac{u}{p}\Big)=\Big( \dfrac{4u}{p}\Big)=\Big( \dfrac{4u+u^2+4}{p}\Big)=\Big( \dfrac{(u+2)^2}{p}\Big)=1 \eqno{\ensuremath{\square}}$$
\vk.3cm
{\bf Lemma\ 4.2.}\quad {\it If $(3.5)$ is consistent and $u_1\equiv \pm 1\ (\mathrm{mod}\ 8)$, then $\Big( \dfrac{n}{p}\Big)=\Big( \dfrac{u_2}{p}\Big)$ for any prime $p$.}
\vk.2cm
{\it Proof}.\quad Indeed, in this case by Lemma 3.5 $\nu_2(z)>\nu_2(x)$. Hence $\nu_2(z-x)=\nu_2(x)$, so we have in $(3.5')$\ $n^m=u_2^k$ with $k, m$ odd, and therefore the conclusion of Lemma 4.2. \hfill\ensuremath{\square}
\vk.3cm
We are ready now to prove Theorem 1.3. Let $p\mid (u^2+4)$. By taking $\Big( \dfrac{}{p}\Big)$ on $(3.5)$ and using Lemmas 4.1, 4.2 one sees that
  $$\Big( \dfrac{u_1}{p}\Big)^x=\Big( \dfrac{n}{p}\Big)^{y-z}=\Big( \dfrac{u_2}{p}\Big)^{y-z}=\Big( \dfrac{u_2}{p}\Big)^y\leqno{(4.1)}$$
(as $z$ is even). Now taking the product of $(4.1)$ over all (not necessarily distinct) prime divisors $p\mid (u^2+4)$ we have
  $$\Big( \dfrac{u_1}{u^2+4}\Big)^x=\prod_{p\mid (u^2+4)}\ \Big( \dfrac{u_1}{p}\Big)^x=\prod_{p\mid (u^2+4)}\ \Big( \dfrac{u_2}{p}\Big)^y=\Big( \dfrac{u_2}{u^2+4}\Big)^y\leqno{(4.2)}$$

By the quadratic reciprocity law
  $$\Big( \dfrac{u_1}{u^2+4}\Big)=\Big( \dfrac{u^2+4}{u_1}\Big)=\Big( \dfrac{2}{u_1}\Big)=1\leqno{(4.3)}$$
  $$\Big( \dfrac{u_2}{u^2+4}\Big)=\Big( \dfrac{u^2+4}{u_2}\Big)=\Big( \dfrac{2}{u_2}\Big)=-1\leqno{(4.4)}$$
as $u_1\equiv \pm 1\ (\mathrm{mod}\ 8)$,\ $u_2\equiv \pm 5\ (\mathrm{mod}\ 8)$. Altogether $(4.2)-(4.4)$ imply that $\Big( \dfrac{u_2}{p}\Big)^y=(-1)^y=1$, {\it i.e.} 
 $y$ must be even. \hfill\ensuremath{\square} 
\vk.3cm
{\bf Corollary\ 4.3.}\quad {\it The possibility $u_1\equiv 7\ (\mathrm{mod}\ 8)$ in Lemma 3.5 is not realized.}
\vk.2cm
{\it Proof}.\quad Indeed, in this case $\nu_2(z)>\nu_2(x)\geq 1$, so $(3.5)$ is inconsistent by Lemma 3.6.\hfill\ensuremath{\square}
\vk.3cm
{\bf Corollary\ 4.4.}\quad {\it In the case $u_1\equiv 1\ (\mathrm{mod}\ 8)$ of Lemma 3.5 we have 
  $$(\nu_2(x), \nu_2(y), \nu_2(z))=(0, \geq 2, 1).$$}

{\it Proof}.\quad By Lemma 3.5 and Theorem 1.3:\ $y, z$ are even, hence $x$ is odd by Lemma 3.6. From the proof of Corollary 3.8 it follows that $\nu_2(z)=1$. For a prime $p\mid (u^2+4)$ by taking $\Big( \dfrac{}{p}\Big)$ on $A=t_1^x$ of 
$(3.14)$ and using Lemma 4.1 one gets 
  $$\Big( \dfrac{t_1}{p}\Big)=\Big( \dfrac{n}{p}\Big)^{(y-z)/2}\leqno{(4.5)}$$

By the same reason of $(4.4)$ we have $\Big( \dfrac{t_1}{u^2+4}\Big)=-1$, as $t_1\equiv 5\ (\mathrm{mod}\ 8)$ by Corollary 3.8. Hence there exists a prime $p_0\mid (u^2+4)$ such that
  $$\Big( \dfrac{t_1}{p_0}\Big)=-1\leqno{(4.6)}$$

From $(4.5), (4.6)$ one concludes that $(y-z)/2$ must be odd (and $\Big( \dfrac{n}{p_0}\Big)=-1$), so the conclusion of Corollary 4.4 follows.
\hfill\ensuremath{\square}  
\vk.3cm
{\bf Remark\ 4.5.}\quad One can have another proof of Lemma 3.6 as shown in several steps below. Assuming $y, z$ even, and arguing as in the proof of Corollary 3.7 one gets equations $(3.14)$ together with $(3.11)-(3.13)$. 
\vk.2cm
1)\ If $u_1\equiv 5\ (\mathrm{mod}\ 8)$ we have four possibilities for $(t_1,t_2)$:
  $$(i)\ t_1\equiv 1\ (\mathrm{mod}\ 8),\ t_2\equiv 5\ (\mathrm{mod}\ 8);\ \ (ii)\ t_1\equiv 5\ (\mathrm{mod}\ 8),\ t_2\equiv 1\ (\mathrm{mod}\ 8);$$
  $$(iii)\ t_1\equiv 3\ (\mathrm{mod}\ 8),\ t_2\equiv 7\ (\mathrm{mod}\ 8);\ \ (iv)\ t_1\equiv 7\ (\mathrm{mod}\ 8),\ t_2\equiv 3\ (\mathrm{mod}\ 8);$$
all of them violate $(3.13)$.
\vk.2cm
2)\ Assume now $u_1\equiv \pm 1\ (\mathrm{mod}\ 8)$ and $x$ {\it even}, hence $\nu_2(z)\geq 2$ by Lemma 3.5. We will shows that $\nu_2(y)=1$. Indeed, considering $p\mid (u^2+4)$ and taking $\Big( \dfrac{}{p}\Big)_4$ on $(3.5)$ 
one has by using Lemmas 4.1, 4.2
  $$\Big( \dfrac{u_1}{p}\Big)^{x/2}=\Big( \dfrac{-1}{p}\Big)_4 \Big( \dfrac{n}{p}\Big)^{(y-z)/2}=\begin{cases}\quad \Big( \dfrac{u_2}{p}\Big)^{y/2},\ p\equiv 1\ (\mathrm{mod}\ 8)\\
      -\Big( \dfrac{u_2}{p}\Big)^{y/2},\ p\equiv 5\ (\mathrm{mod}\ 8)\end{cases}\leqno{(4.7)} $$
as $z/2$ is even. Let $r$ denote the number of prime divisors $p\mid (u^2+4)$,\ $p\equiv 5\ (\mathrm{mod}\ 8)$. Clearly $r$ is {\it odd}, as $u^2+4\equiv 5\ (\mathrm{mod}\ 8)$. 
In a similar way as in $(4.2)-(4.4)$, taking the product of $(4.7)$ over all (not necessarily distinct) prime divisors $p\mid (u^2+4)$ we get
  $$1=\Big( \dfrac{u_1}{u^2+4}\Big)^{x/2}=(-1)^r \Big( \dfrac{u_2}{u^2+4}\Big)^{y/2}=-(-1)^{y/2}.$$

Hence $y/2$ must be odd, so $(y-z)/2$ is odd. For any prime $p\mid (u^2+4)$ taking $\Big( \dfrac{}{p}\Big)$ on equation $A=t_1^x$ from $(3.14)$ now gives us
  $$\Big( \dfrac{n}{p}\Big)=1\ \bigg(=\Big( \dfrac{u_2}{p}\Big)\ \ \mathrm{by\ Lemma\ 4.2}\bigg)\leqno{(4.8)}$$

On the other hand from $(4.4)$ it follows that there exists a prime $p_0\mid (u^2+4)$ such that $\Big( \dfrac{u_2}{p_0}\Big)=-1$, a contradiction with $(4.8)$. Thus $(3.14)$ (and hence $(3.5)$) is inconsistent.   

\vk.5cm

\begin{center}
{\bf 5.\ The case $u_1\equiv 5\ (\mathrm{mod}\ 8)$}
\end{center}

In this case by $(3)$ of Lemma 3.5 we have $\nu_2(z)=\nu_2(x)$, hence from $(3.4)$ it follows that $u_2=w^{2^s}$, where $s=\nu_2(z-x)-\nu_2(x)$. The following lemma can be proved similarly as Lemma 4.2.
\vk.3cm
{\bf Lemma\ 5.1.}\quad {\it If $(3.5)$ is consistent and $u_1\equiv 5\ (\mathrm{mod}\ 8)$, then $\Big( \dfrac{n}{p}\Big)=\Big( \dfrac{w}{p}\Big)$ for any prime $p$.}
\vk.2cm
{\it Proof}.\quad Indeed, in this case $n^m=w^k$ with $k, m$ odd by the above argument, and therefore the conclusion of Lemma 5.1. \hfill\ensuremath{\square}

\vk.3cm
{\bf Lemma\ 5.2.}\quad {\it If $x, z$ are even and $(3.5)$ is consistent, then $y$ is odd and $u_1\equiv 5\ (\mathrm{mod}\ 8)$. Moreover 
$n$ admits a decomposition $n=n_{1} n_{2}$ such that $\gcd(n_{1},n_{2})=1$ and 
  $$\begin{cases} &\qquad u_1^{x/2}\quad =u^y n_{2}^{y-z}-2^{2y-2} n_{1}^{y-z}\\
      & (u^2+4)^{z/2}=u^y n_{2}^{y-z}+2^{2y-2} n_{1}^{y-z}\end{cases}\leqno{(5.1)}$$
}

{\it Proof}.\quad By Lemma 3.6 $y$ is odd. In view of Lemma 3.5 and Theorem 1.3 we are in the situation $(3)$ of Lemma 3.5. Now one rewrites $(3.5)$ as
  $$C_1\cdot D_1=2^{2y} u^y n^{y-z}\ \mathrm{with}\ \gcd(C_1,D_1)=2,\  2\| D_1,$$
where
  $$C_1=(u^2+4)^{z/2}-u_1^{x/2},\  D_1=(u^2+4)^{z/2}+u_1^{x/2}.$$

As $2\| D_1$ we obtain either 
  $$C_1=2^{2y-1} n_{1}^{y-z},\ D_1=2 u^y n_{2}^{y-z}\leqno{(5.2)}$$
or 
  $$C_1=2^{2y-1} u^y n_{1}^{y-z},\ D_1=2 n_{2}^{y-z}\leqno{(5.3)}$$
where $n=n_{1} n_{2}$,\ $\gcd(n_{1},n_{2})=1$ and 
  $$w=w_{1} w_{2},\ n_{1}^m=w_{1}^k,\ n_{2}^m=w_{2}^k\leqno{(5.4)}$$
with $k, m$ odd from Lemma 5.1. Note that this is not used in the proof here, we label it for convenience in proving Proposition 5.5 below. 
\vk.2cm
Clearly $(5.2)$ is equivalent to $(5.1)$. It remains to show that $(5.3)$ can't happen by rewriting it as
  $$\begin{cases} &\qquad u_1^{x/2}\quad =n_{2}^{y-z}-2^{2y-2} u^y n_{1}^{y-z}\\
      & (u^2+4)^{z/2}=n_{2}^{y-z}+2^{2y-2} u^y n_{1}^{y-z}\end{cases}\leqno{(5.5)}$$
which is impossible, since $(u^2+4)^{z/2} < 2^{2y-2} u^y$. 
\hfill\ensuremath{\square}
\vk.3cm
{\bf Lemma\ 5.3.}\quad {\it If $\Big( \dfrac{u_1}{u}\Big)=1$ and $u_2$ is a square, then $u\equiv 1\ (\mathrm{mod}\ 4)$.}
\vk.2cm
{\it Proof}.\quad We have obviously
  $$1=\Big( \dfrac{u_1}{u}\Big)=\Big( \dfrac{u_1 u_2}{u}\Big)=\Big( \dfrac{u^2-4}{u}\Big)=\Big( \dfrac{-1}{u}\Big),$$
so the conclusion of the lemma. \hfill\ensuremath{\square}
\vk.3cm
{\bf Lemma\ 5.4.}\quad {\it In the notations of Lemma 5.1 we have
\vk.1cm
$(1)$\ if $w\equiv \pm 3\ (\mathrm{mod}\ 8)$, then $x, z$ are odd, $y$ is even$;$
\vk.1cm
$(2)$\ if $w\equiv \pm 1\ (\mathrm{mod}\ 8)$, then $x, z$ are even, $y$ is odd.}
\vk.2cm
{\it Proof}.\quad For a prime $p\mid (u^2+4)$ by taking $\Big( \dfrac{}{p}\Big)$ on $(3.5)$ and using Lemmas 4.1, 5.1 one sees that
  $$\Big( \dfrac{u_1}{p}\Big)^x=\Big( \dfrac{n}{p}\Big)^{y-z}=\Big( \dfrac{w}{p}\Big)^{y-z}\leqno{(5.6)}$$

By taking the product of both sides of $(5.6)$ over all (not necessarily distinct) prime divisors $p\mid (u^2+4)$ and using the reciprocity law we have
  $$\prod_{p\mid (u^2+4)}\ \Big( \dfrac{u_1}{p}\Big)^x=\Big( \dfrac{u_1}{u^2+4}\Big)^x=\Big( \dfrac{u^2+4}{u_1}\Big)^x=\Big( \dfrac{2}{u_1}\Big)^x=(-1)^x\leqno{(5.7)}$$
  $$\begin{aligned}
\prod_{p\mid (u^2+4)}\ \Big( \dfrac{w}{p}\Big)^{y-z}& =\Big( \dfrac{w}{u^2+4}\Big)^{y-z}=\Big( \dfrac{u^2+4}{w}\Big)^{y-z}=\\
& =\Big( \dfrac{2}{w}\Big)^{y-z}=\begin{cases} & (-1)^{y-z},\ w\equiv \pm 3
\ (\mathrm{mod}\ 8)\\
&\quad 1,\qquad \ w\equiv \pm 1\ (\mathrm{mod}\ 8)\end{cases}
\end{aligned}
\leqno{(5.8)}$$

Hence if $w\equiv \pm 3\ (\mathrm{mod}\ 8)$, then by equalizing $(5.7), (5.8)$:\ $(-1)^x=(-1)^{y-z}$. Thus $y$ must be even, as $\nu_2(z)=\nu_2(x)$. In view of Lemma 3.6 $x, z$ are odd.
\vk.2cm
In the case $w\equiv \pm 1\ (\mathrm{mod}\ 8)$, again equalizing $(5.7), (5.8)$ we see that \ $(-1)^x=1$, therefore $x$ is even, and so is $z$. By Lemma 3.6 $y$ must be odd. \hfill\ensuremath{\square}
\vk.3cm
{\bf Proposition\ 5.5.}\quad {\it In the situation of Lemma 5.4 we have
\vk.1cm
$(1)$\ if $w\equiv \pm 3\ (\mathrm{mod}\ 8)$, then $u\equiv 1\ (\mathrm{mod}\ 4);$
\vk.1cm
$(2)$\ if $w\equiv \pm 1\ (\mathrm{mod}\ 8)$, then $u\equiv\pm 3\ (\mathrm{mod}\ 8)$. Moreover, if $u\equiv 3\ (\mathrm{mod}\ 8)$,\ then $w$ can not be a square.}
\vk.2cm
{\it Proof}.\quad $(1)$\ If $w\equiv \pm 3\ (\mathrm{mod}\ 8)$, then $x, z$ are odd in view of Lemma 5.4. So by taking $\Big( \dfrac{}{u}\Big)$ on $(3.5)$ one gets $\Big( \dfrac{u_1}{u}\Big)=1$, hence 
$u\equiv 1\ (\mathrm{mod}\ 4)$ by Lemma 5.3.
\vk.2cm
$(2)$\ In the case $w\equiv \pm 1\ (\mathrm{mod}\ 8)$:\ $x, z$ are even, $y$ is odd by Lemma 5.4. There are two subcases to consider.
\vk.3cm
$\boldsymbol{I.}$\ {\it $x/2, z/2$ are odd}.\quad For a prime $p\mid (u^2+4)$ by taking $\Big( \dfrac{}{p}\Big)$ on $D_1=2 u^y n_{2}^{y-z}$ from $(5.2), (5.4)$ and using Lemmas 4.1, 5.1 one sees that
  $$\Big( \dfrac{u_1}{p}\Big)=\Big( \dfrac{2}{p}\Big) \Big( \dfrac{n_{2}}{p}\Big)=\Big( \dfrac{2}{p}\Big) \Big( \dfrac{w_{2}}{p}\Big)=\begin{cases} &\ \ \Big( \dfrac{w_{2}}{p}\Big),\ \ p\equiv 1
\ (\mathrm{mod}\ 8)\\
& -\Big( \dfrac{w_{2}}{p}\Big),\ \ p\equiv 5\ (\mathrm{mod}\ 8)\end{cases}\leqno{(5.9)}$$

Recall that the number of (not necessarily distinct) prime divisors $p\mid (u^2+4)$,\ $p\equiv 5\ (\mathrm{mod}\ 8)$ is odd, so $\und{p\mid (u^2+4)}\prod\ \Big( \dfrac{2}{p}\Big)=-1$. Now taking the product of both sides of $(5.9)$ 
over all (not necessarily distinct) prime divisors $p\mid (u^2+4)$ and using the reciprocity law one has
  $$\prod_{p\mid (u^2+4)}\ \Big( \dfrac{u_1}{p}\Big)=\Big( \dfrac{u_1}{u^2+4}\Big)=\Big( \dfrac{u^2+4}{u_1}\Big)=\Big( \dfrac{2}{u_1}\Big)=-1\leqno{(5.10)}$$
and
  $$\prod_{p\mid (u^2+4)}\ \Big( \dfrac{2}{p}\Big) \Big( \dfrac{w_{2}}{p}\Big)=-\prod_{p\mid (u^2+4)}\ \Big( \dfrac{w_{2}}{p}\Big)=-\Big( \dfrac{w_{2}}{u^2+4}\Big)=-\Big( \dfrac{u^2+4}{w_{2}}\Big)=-\Big( \dfrac{2}{w_{2}}\Big)\leqno{(5.11)}$$

Equalizing $(5.10)$, $(5.11)$ we get $w_{2}\equiv\pm 1\ (\mathrm{mod}\ 8)$, so in view of $(5.4)$:\ $n_{2}\equiv\pm 1\ (\mathrm{mod}\ 8)$. From this and $(5.1)$ it follows that $u\equiv\pm 3\ (\mathrm{mod}\ 8)$. Moreover, if $u\equiv 3\ (\mathrm{mod}\ 8)$,\ then 
$w_{2}\equiv -1\ (\mathrm{mod}\ 8)$, hence by $(5.4)$\ $w$ can not be a square.
\vk.3cm
$\boldsymbol{II.}$\ {\it $x/2, z/2$ are even}.\quad If one takes $\Big( \dfrac{}{u}\Big)$ on the second equation of $(5.1)$, then $\Big( \dfrac{n_{1}}{u}\Big)=1$. Now taking $\Big( \dfrac{}{u}\Big)$ on the first 
equation of $(5.1)$ we get $1=\Big( \dfrac{-1}{u}\Big) \Big( \dfrac{n_{1}}{u}\Big)$. Thus $u\equiv 1\ (\mathrm{mod}\ 4)$.
\vk.2cm
The proof of Proposition 5.5 is completed.  \hfill\ensuremath{\square}
\vk.3cm
As for Theorem 1.4  notice that the case $u_1\equiv\pm 1\ (\mathrm{mod}\ 8)$ follows from Corollaries 3.7, 3.8, 4.3 and 4.4. The rest of Theorem 1.4, {\it i.e.} the case $u_1\equiv 5\ (\mathrm{mod}\ 8)$, follows from 
Lemma 5.4 and Proposition 5.5. 
\vk.2cm
The equalities for Jacobi symbols are immediate from $(3.5)$ and Lemma 5.1. \hfill\ensuremath{\square}

\begin{center}
{\bf 6.\ Proof of Corollary 1.5}
\end{center}

In this section we shall apply results of previous parts for establishing the truth of Je{\'s}manowicz' conjecture for $u<100$ and $v = 2$. In view of Theorem 1.4 one has to consider only two cases:\ $u_1\equiv 1\ (\mathrm{mod}\ 8)$ 
and $u_1\equiv 5\ (\mathrm{mod}\ 8)$.
\vk.3cm
{\bf Observation\ 6.1.}\quad {\it If $u_1\equiv 1\ (\mathrm{mod}\ 8)$ and $(3.5)$ is consistent, then $u>183$.}
\vk.2cm
{\it Proof}.\quad Indeed, it was noted that $x\geq 3$ by $(3.14)$. On the other hand from the proof of Corollary 3.8 we have $\nu_2(z)=1$, so $z\geq 6$, hence $y\geq 8$. From $(3.12)$ it follows that $2^{y+1}\mid t_2-t_1$, 
as $x$ is odd. Since $t_1, u_2$ are co-prime and $\equiv 5\ (\mathrm{mod}\ 8)$, so $t_1 u_2\geq 5\cdot 13$. Therefore $u>\sqrt{t_1 t_2 u_2}\geq \sqrt{(2^9+5)\cdot 65}>183$. \hfill\ensuremath{\square}
\vk.3cm
{\bf Observation\ 6.2.}\quad {\it If $u_1\equiv 1\ (\mathrm{mod}\ 8)$ and $(3.5)$ is consistent, then in fact $u>729$.}
\vk.2cm
{\it Proof}.\quad By Corollary 4.4 one knows $4\mid y$. We claim that $y\geq 12$. Assuming on the contrary $y=8$, then by the above $z=6$. In view of $(3.11)$ and \cite{Wa58} we must have $x>3$, so $x=5$, which gives us a non-trivial 
solution of $X^5+Y^5=2 Z^3$. This is impossible by \cite{BeVY04} (Theorem 1.5). 
\vk.2cm
Therefore $y\geq 12$, and by the argument above $u>\sqrt{(2^{13}+5)\cdot 65}>729$.
\hfill\ensuremath{\square}
\vk.2cm
It remains to consider the case $u_1\equiv 5\ (\mathrm{mod}\ 8)$. In the range of odd primes $<100$ there are ten possibilities with $u^2-4=u_1 u_2$ and $u_2$ is a square, namely $u=7, 11, 23$, $43, 47, 61, 73, 79, 83, 97.$ 
In view of Proposition 5.5 we shall exclude the possibilities $u=7, 23, 47, 79$.
\vk.3cm
{\bf 6.3.}\ For $(u,u_1,u_2)=(11,13,3^2)$,\ $(43,5\cdot 41,3^2)$,\ $(83,5\cdot 17,3^4)$  we have $w\equiv \pm 3\ (\mathrm{mod}\ 8)$, hence $u\equiv 1\ (\mathrm{mod}\ 4)$ by 
Proposition 5.5, a contradiction. Note that in the original version to eliminate the possibility $(83,5\cdot 17,3^4)$ and $w=9$ we used implicitly the fact that if $u\equiv 3\ (\mathrm{mod}\ 8)$, then $w$ can not be a square, which we include a proof 
in the revised version ({\it cf.} Proposition 5.5 above). The referee provides another argument by choosing $p=5\mid u_1$ which leads also to a contradiction as follows
    $$1=\Big( \dfrac{9}{5}\Big)=\Big( \dfrac{w}{p}\Big)\ne \Big( \dfrac{u}{p}\Big)=\Big( \dfrac{83}{5}\Big)=-1.$$
\vk.3cm
{\bf 6.4.}\ For $(u,u_1,u_2)=(61,7\cdot 59,3^2)$\ one has $w=3$, so
  $$-1=\Big( \dfrac{w}{7}\Big)\ne \Big( \dfrac{u^2+4}{7}\Big)=1 $$
a contradiction with $(2.1)$ of Theorem 1.4.
\vk.3cm
{\bf 6.5.}\ For $(u,u_1,u_2)=(73,3\cdot 71,5^2)$\ we have $w=5$, hence $x, z$ are odd and $y$ is even by $(2.1)$ of Theorem 1.4. Taking modulo $73$ on $(3.5)$ one gets 
  $$4^z\equiv (-6)^x\ (\mathrm{mod}\ 73)\leqno{(6.1)}$$

Working in $\F_{73}^*$ we have
  $$\mathrm{ord}(4)=9,\quad \mathrm{ord}(-6)=36\leqno{(6.2)}$$

Therefore from $(6.1), (6.2)$ it follows that $36\mid 9 x$, so $4\mid x$, a contradiction. 
\vk.3cm
{\bf 6.6.}\ For $(u,u_1,u_2)=(97,5\cdot 11\cdot 19,3^2)$\ one has $w=3$, so
    $$1=\Big( \dfrac{w}{11}\Big)\ne \Big( \dfrac{u^2+4}{11}\Big)=-1 $$
again a contradiction with $(2.1)$ of Theorem 1.4.
\vk.5cm

\vk.5cm
FPT University, Department of Mathematics, Hoa Lac High Tech Park, Hanoi, Vietnam
\vk.2cm
Email addresses:\ thiennv15@fe.edu.vn
\vk.3cm
FPT University, Department of Mathematics and Information Assurance, Hoa Lac High Tech Park, Hanoi, Vietnam
\vk.2cm
Email addresses:\ vietnk@fe.edu.vn
\vk.3cm
FPT University, Department of Mathematics and Office of Science Management and International Affairs, Hoa Lac High Tech Park, Hanoi, Vietnam
\vk.2cm
Email addresses:\ quyph@fe.edu.vn

\end{document}